\documentclass[10pt,twoside]{article}
\def\YEAR{\year}\newcount\VOL\VOL=\YEAR\advance\VOL by-1995
\def\firstpage{1}\def\lastpage{1000}
\def\received{}\def\revised{}
\def\communicated{}

\makeatletter
\def\magnification{\afterassignment\m@g\count@}
\def\m@g{\mag=\count@\hsize6.5truein\vsize8.9truein\dimen\footins8truein}
\makeatother

\font\eightrm=cmr8
\font\caps=cmcsc10                    
\font\Caps=cmcsc10 scaled \magstep1   

\def\bfseries{\normalsize\caps}

\makeatletter
\@specialpagefalse\headheight=8.5pt
\def\DocMath{}
\renewcommand{\@evenhead}{%
    \ifnum\thepage>\lastpage\rlap{\thepage}\hfill%
    \else\rlap{\thepage}\slshape\leftmark\hfill{\caps\SAuthor}\hfill\fi}%
\renewcommand{\@oddhead}{%
    \ifnum\thepage=\firstpage{\DocMath\hfill\llap{\thepage}}%
    \else{\slshape\rightmark}\hfill{\caps\STitle}\hfill\llap{\thepage}\fi}%
\makeatother

\def\TSkip{\bigskip}
\newbox\TheTitle{\obeylines\gdef\GetTitle #1
\ShortTitle  #2
\SubTitle    #3
\Author      #4
\ShortAuthor #5
\EndTitle
{\setbox\TheTitle=\vbox{\baselineskip=20pt\let\par=\cr\obeylines%
\halign{\centerline{\Caps##}\cr\noalign{\medskip}\cr#1\cr}}%
        \copy\TheTitle\TSkip\TSkip%
\def\next{#2}\ifx\next\empty\gdef\STitle{#1}\else\gdef\STitle{#2}\fi%
\def\next{#3}\ifx\next\empty%
    \else\setbox\TheTitle=\vbox{\baselineskip=20pt\let\par=\cr\obeylines%
    \halign{\centerline{\caps##} #3\cr}}\copy\TheTitle\TSkip\TSkip\fi%
\centerline{\caps #4}\TSkip\TSkip%
\def\next{#5}\ifx\next\empty\gdef\SAuthor{#4}\else\gdef\SAuthor{#5}\fi%
\ifx\received\empty\relax
    \else\centerline{\eightrm Received: \received}\fi%
\ifx\revised\empty\TSkip%
    \else\centerline{\eightrm Revised: \revised}\TSkip\fi%
\ifx\communicated\empty\relax
    \else\centerline{\eightrm Communicated by \communicated}\fi\TSkip\TSkip%
\catcode'015=5}}\def\Title{\obeylines\GetTitle}
\def\Abstract{\begingroup\narrower
    \parskip=\medskipamount\parindent=0pt{\caps Abstract. }}
\def\EndAbstract{\par\endgroup\TSkip}

\long\def\MSC#1\EndMSC{\def\arg{#1}\ifx\arg\empty\relax\else
     {\par\narrower\noindent%
     2000 Mathematics Subject Classification: #1\par}\fi}

\long\def\KEY#1\EndKEY{\def\arg{#1}\ifx\arg\empty\relax\else
        {\par\narrower\noindent Keywords and Phrases: #1\par}\fi\TSkip}

\newbox\TheAdd\def\Addresses{\vfill\copy\TheAdd\vfill
    \ifodd\number\lastpage\vfill\eject\phantom{.}\vfill\eject\fi}
{\obeylines\gdef\GetAddress #1
\Address #2 
\Address #3
\Address #4
\EndAddress
{\def\xs{4.3truecm}\parindent=0pt
\setbox0=\vtop{{\obeylines\hsize=\xs#1\par}}\def\next{#2}
\ifx\next\empty 
     \setbox\TheAdd=\hbox to\hsize{\hfill\copy0\hfill}
\else\setbox1=\vtop{{\obeylines\hsize=\xs#2\par}}\def\next{#3}
\ifx\next\empty 
     \setbox\TheAdd=\hbox to\hsize{\hfill\copy0\hfill\copy1\hfill}
\else\setbox2=\vtop{{\obeylines\hsize=\xs#3\par}}\def\next{#4}
\ifx\next\empty\ 
     \setbox\TheAdd=\vtop{\hbox to\hsize{\hfill\copy0\hfill\copy1\hfill}
                \vskip20pt\hbox to\hsize{\hfill\copy2\hfill}}
\else\setbox3=\vtop{{\obeylines\hsize=\xs#4\par}}
     \setbox\TheAdd=\vtop{\hbox to\hsize{\hfill\copy0\hfill\copy1\hfill}
                \vskip20pt\hbox to\hsize{\hfill\copy2\hfill\copy3\hfill}}
\fi\fi\fi\catcode'015=5}}\gdef\Address{\obeylines\GetAddress}

\usepackage[all]{xy}
\usepackage{supertabular}
\usepackage{tom}

\def\theoremnumberstyle{section} 
\def\proofname{Proof}

\renewenvironment{equationlist}{
  \begin{list}{}
    {\renewcommand{\makelabel}[1]{\stepcounter{equation}\mbox{\rm(\theequation)}%
        \if##1\empty\else\mylabel{##1}{\theequation}\fi}
     \leftmargin2cm \itemindent0cm\labelwidth1.5cm \parsep6pt \topsep6pt\labelsep0.2cm}
    }
  {\end{list}
  }
\newenvironment{mylist}{
  \begin{list}{}
    {\leftmargin4ex \itemindent0ex\labelwidth2ex \parsep0.5ex \topsep0pt\labelsep1ex}
    }
  {\end{list}
  }
\newenvironment{myenumerate}{
  \begin{list}{}
    {\leftmargin5ex \itemindent0ex\labelwidth4ex \parsep0.5ex \topsep0pt\labelsep1ex}
    }
  {\end{list}
  }

\numberwithin{equation}{section} 
\allowdisplaybreaks[3]   
\newcommand{\asf}{1}                                
\renewcommand{\arraystretch}{\asf}                    
\newcommand{\as}[1]{\renewcommand{\arraystretch}{#1}} 
\newcommand{\expl}[1]{{_{_{\mbox{\tiny #1}}}}}  
\SelectTips{cm}{12}  
\CompileMatrices     

\begin{document}
\Title 
Irreducibility of Equisingular Families of Curves 
Improved Conditions
\ShortTitle 
Irreducibility
\SubTitle  
\Author 
Thomas Keilen
\ShortAuthor 
Thomas Keilen
\EndTitle
\Abstract 
In \cite{Kei03} we gave sufficient conditions for the
irreducibility of the family
$V_{|D|}^{irr}\big(\ks_1,\ldots,\ks_r\big)$ of irreducible curves
in the linear system $|D|_l$ with precisely $r$ singular points of
topological respectively analytical types $\ks_1,\ldots,\ks_r$ on
several classes of smooth projective surfaces $\Sigma$. The conditions
where of the form
\begin{displaymath}
  \sum\limits_{i=1}^r\big(\tau^*(\ks_i)+2\big)^2
  <
  \gamma\cdot (D- K_\Sigma)^2,       
\end{displaymath}
where $\tau^*$ is some invariant of singularity types,
$K_\Sigma$ is the canonical divisor of $\Sigma$ and $\gamma$ is
some constant. In the present paper we improve this condition,
that is the constant $\gamma$, by a factor $9$. 
\EndAbstract

\MSC 
Primary 14H10, 14H15, 14H20; Secondary 14J26, 14J27, 14J28, 14J70
\EndMSC
\KEY 
Algebraic geometry, singularity theory
\EndKEY
\Address
Thomas Keilen
Mathematics Institute
University of Warwick
Coventry CV4 7AL

\Address
\Address
\Address

\EndAddress




   \tableofcontents


   \section{Introduction}

   If we fix a linear system $|D|_l$ on a smooth
   projective surface $\Sigma$ over $\C$ and singularity types
   $\ks_1,\ldots,\ks_r$ we denote by
   $V^{irr}=V_{|D|}^{irr}\big(\ks_1,\ldots,\ks_r\big)$ the variety 
   of irreducible curves in $|D|_l$ with
   precisely $r$ singular points of the given types. We would like to
   give numerical conditions, depending on the singularity types, the
   linear system and the surface, which ensure that the family
   $V$ is irreducible,
   once it is non-empty.

   In order to keep the presentation as short as possible we refer the
   reader to \cite{Kei03} for an introduction to the significance of
   the question and for most of the notation we are going to
   use. Moreover, we will apply many of the technical results shown
   there. The proof runs along the same lines as the original one by
   showing that some irreducible ``regular'' subscheme of $V$ is dense
   in $V$. We do this again by considering a morphism $\Phi$ on a certain
   subscheme of $V$ and comparing dimensions. However, the subscheme
   which we consider and the morphism are completely different. 

   We now
   introduce these new objects. In Section \ref{sec:irred} we then
   formulate the main results, and we prove them in Section
   \ref{sec:technical-lemmata}. Lemma \ref{lem:irred-E-fix} is the most
   important technical adjustment which leads to the improved
   coefficient.



     \subsection{The Deformation Determinacy}

     If $\ks$ is a topological (respectively analytical) singularity
     type with representative $(C,z)$ then 
     \begin{displaymath}
       \nu^s(\ks)=\nu^s(C,z)=\min\big\{m\geq
       0\;\big|\;\m_{\Sigma,z}^{m+1}\subseteq I^s(C,z)\big\}
     \end{displaymath}
     respectively
     \begin{displaymath}
       \nu^a(\ks)=\nu^a(C,z)=\min\big\{m\geq
       0\;\big|\;\m_{\Sigma,z}^{m+1}\subseteq I^a(C,z)\big\},
     \end{displaymath}
     where $I^s(C,z)=\kj_{X^s(C)/\Sigma,z}$ is the singularity ideal
     of the topological singularity type $(C,z)$ and $I^a(C,z)=\kj_{X^a(C)/\Sigma,z}$ is the
     analytical singularity ideal of $(C,z)$ respectively
     (cf. \cite{Kei03} Section 1.3). These are invariants of the
     topological (respectively analytical) singularity type satisfying (cf. \cite{GLS00} Section 1.2 and 1.3)
     \begin{displaymath}
       \nu^s(\ks)\leq \tau^{es}(\ks)\;\;\;\;\mbox{ respectively
         }\;\;\;\;
       \nu^a(\ks)\leq \tau(\ks),
     \end{displaymath}
     and they
     are called \emph{topological deformation determinacy}
     (respectively \emph{analytical deformation
       determinacy})

   \subsection{Singularity Schemes}\label{subsec:schemes}          
     For a reduced curve $C\subset\Sigma$ we recall the definition of
     the zero-dimensional schemes $X^{es}_{fix}(C)$ and
     $X^{ea}_{fix}(C)$ from \cite{GLS00} Section 1.1. They are defined
     by the ideal sheaves $\kj_{X^{es}_{fix}(C)/\Sigma}$
     and $\kj_{X^{ea}_{fix}(C)/\Sigma}$
     respectively, given by the following stalks 
     \begin{mylist}
     \item[$\bullet$] $\kj_{X^{es}_{fix}(C)/\Sigma,z}=I^{es}_{fix}(C,z)=
       \left\{g\in\ko_{\Sigma,z}\;\Big|\;
         {\tiny
         \begin{array}{c}
           f+\varepsilon g \mbox{ is equisingular
             over } \C[\varepsilon]/(\varepsilon^2)\\
           \mbox{ along the
             trivial section}
         \end{array}
         }
       \right\}$,
       where $f\in\ko_{\Sigma,z}$ is a local equation of $C$ at $z$.
     \item[$\bullet$] $\kj_{X^{ea}_{fix}(C)/\Sigma,z}=I^{ea}_{fix}(C,z)=
       \langle f\rangle +\m\cdot\big\langle \tfrac{\partial f}{\partial x},\tfrac{\partial
         f}{\partial y}\big\rangle\subseteq\ko_{\Sigma,z}$, 
       where $x,y$ denote local coordinates of $\Sigma$ at $z$ and
       $f\in\ko_{\Sigma,z}$ is a local equation of $C$.
     \end{mylist}
     So by definition we have 
     \begin{displaymath}
       \deg\big(X^{es}_{fix}(C),z\big)=\tau^{es}(C,z)+2
       \;\;\;\;\mbox{ and }\;\;\;\;
       \deg\big(X^{ea}_{fix}(C),z\big)=\tau(C,z)+2.
     \end{displaymath}
     \begin{center}
       \framebox[11cm]{
         \begin{minipage}{10cm}
         \medskip
           Throughout this article we will frequently treat
           topological and analytical singularities at the same time.
           Whenever we do so, we will write $X^*_{fix}(C)$ for $X^{es}_{fix}(C)$
           respectively for $X^{ea}_{fix}(C)$, we will write
           $\nu^*(\ks)$ for $\nu^s(\ks)$ respectively
           $\nu^a(\ks)$, and we will write $\tau^*(\ks)$ for
           $\tau^{es}(\ks)$ respectively $\tau(\ks)$. For the schemes
           borrowed from 
           \cite{Kei03} we stick to the analogous convention
           made there.
         \medskip
         \end{minipage}
         }
     \end{center}

   \subsection{Equisingular Families and Fibrations}\label{subsec:families}
     Given a divisor $D\in\Div(\Sigma)$ and topological (respectively analytical) singularity types
     $\ks_1,\ldots, \ks_r$. 

     We denote
     by 
     $V^{irr,fix}=V_{|D|}^{irr,fix}(\ks_1,\ldots,\ks_r)$ the
     open  subscheme  of $V^{irr}$ given as
     \begin{displaymath}
       V^{irr,fix}=
       \big\{C\in V^{irr}_{|D|}(\ks_1,\ldots,\ks_r)
       \;\big|\;h^1\big(\Sigma,\kj_{X^*_{fix}(C)/\Sigma}(D)\big)=0\big\}.
     \end{displaymath}

     We define the fibration
     $\Phi=\Phi_D(\ks_1,\ldots,\ks_r)$ by
     \begin{displaymath}
       \xymatrix@R=0.2cm{
         \Phi:
         V_{|D|}^{irr}(\ks_1,\ldots,\ks_r)\ar[r]&
         \Sym^r(\Sigma):
         C\ar@{|->}[r]&
         \;\Sing(C),
         }         
     \end{displaymath}
     sending a curve $C$ to the unordered tuple of its singular
     points.
     
     Note that 
     $H^0\big(\Sigma,\kj_{X^*_{fix}(C)}(D)\big)/H^0(\ko_\Sigma)$ is
     the tangent space of the fibre $\Phi^{-1}\big(\Phi(C)\big)$ has
     at $C\in V^{irr}$, so     
     that 
     \begin{equation}\label{eq:psi-fix}
       \dim\Big(\Phi^{-1}\big(\Phi(C)\big)\Big)\leq
       h^0\big(\Sigma,\kj_{X^*_{fix}(C)}(D)\big) -1.
     \end{equation}
     Moreover, suppose that 
     $h^1\big(\Sigma,\kj_{X^*_{fix}(C)/\Sigma}(D)\big)=0$, then the
     germ of the fibration at $C$ 
     \begin{displaymath}
       (\Phi,C):\big(V,C\big)\rightarrow \big(\Sym^r(\Sigma),\Sing(C)\big)
     \end{displaymath}
     is smooth of fibre dimension
     $h^0\big(\Sigma,\kj_{X^*_{fix}(C)}(D)\big) -1$, i.~e.~ locally
     at $C$ the morphism $\Phi$ is a projection of the product of
     the smooth base space with the smooth fibre.
     This implies in particular, that close to $C$ there is a curve
     having its singularities in very general position.
     (Cf.~\cite{Los98} Proposition 2.1 (e).)


   \section{The Main Results}\label{sec:irred}

   In this section we give sufficient conditions for the
   irreducibility of equisingular families of curves on certain
   surfaces with Picard number one -- including the projective
   plane, general surfaces in $\PC^3$ and general K3-surfaces --, on
   products of curves, and on a subclass of geometrically ruled
   surfaces.

   \subsection{Surfaces with Picard Number One}
   
   \begin{theorem}\label{thm:irred-p3}
     Let $\Sigma$ be a surface such that
     \begin{myenumerate}
     \item[\rm(i)] $\NS(\Sigma)=L\cdot\Z$ with $L$ ample, and
     \item[\rm(ii)] $h^1(\Sigma,C)=0$, whenever $C$ is effective.
     \end{myenumerate}

     Let $D\in\Div(\Sigma)$, let $\ks_1,\ldots,\ks_r$ be 
     topological (respectively  analytical) singularity types.

     Suppose that 
     \begin{equationlist}
     \item[eq:irred-p3:0]  \hspace*{-0.5cm}$D-K_\Sigma$ is big and nef,
     \item[eq:irred-p3:0+]  \hspace*{-0.5cm}$D+K_\Sigma$ is nef,
     \item[eq:irred-p3:1]  \hspace*{-0.5cm}$\sum\limits_{i=1}^r\big(\tau^*(\ks_i)+2\big)<\beta\cdot
       (D-K_\Sigma)^2$\;\;
       for some $0<\beta\leq
       \tfrac{1}{4}$, and 
     \item[eq:irred-p3:2] \hspace*{-0.5cm}$\sum\limits_{i=1}^r\big(\tau^*(\ks_i)+2\big)^2
       <
       \gamma\cdot (D-K_\Sigma)^2$,
       \\\hspace*{4cm}where
       $\gamma=\tfrac{\big(1+\sqrt{1-4\beta}\big)^2\cdot
         L^2}{4\cdot\chi(\ko_\Sigma)+\max\{0,2\cdot K_\Sigma.L\}+6\cdot
         L^2}$.
     \end{equationlist}
     
     Then 
     $V_{|D|}^{irr}(\ks_1,\ldots,\ks_r)$ is empty
     or  irreducible of the expected dimension.
     \hfill $\Box$
   \end{theorem}

   \begin{remark}\label{rem:irred-p3}
     If we set
     \begin{displaymath}
       \gamma=\frac{36\alpha}{(3\alpha+4)^2}\;\;\;\text{ with }\;\;\; 
       \alpha=\frac{4\cdot\chi(\ko_\Sigma)+\max\{0,2\cdot K_\Sigma.L\}+6\cdot
       L^2}{L^2},
     \end{displaymath}
     then a simple calculation shows that
     \eqref{eq:irred-p3:1} becomes redundant.
     For this we have to take into account that
     $\tau^*(\ks)\geq 1$ for any singularity type $\ks$. 
     The claim then follows with
     $\beta=\frac{1}{3}\cdot\gamma\leq \frac{1}{4}$.
     \hfill $\Box$
   \end{remark}

   We now apply the result in several special cases, combining the
   above theorem with the existence results in \cite{KT02} and the
   T-smoothness results in \cite{GLS97}.

   \begin{corollary}
     Let $d\geq 3$, $L\subset\PC^2$ be a line, and 
     $\ks_1,\ldots,\ks_r$ be 
     topological or analytical singularity types.

     Suppose that 
     \begin{displaymath}
       \sum\limits_{i=1}^r \big(\tau^*(\ks_i)+2\big)^2
       <
       \tfrac{90}{289}\cdot (d+3)^2.
     \end{displaymath}

     Then 
     $V_{|dL|}^{irr}(\ks_1,\ldots,\ks_r)$ is non-empty,
     irreducible and T-smooth.
     \hfill $\Box$
   \end{corollary}

   The best general
   results in this case can still be found in \cite{GLS00} (see also \cite{Los98} Corollary
   6.1), where the coefficient on the right hand side is $\frac{9}{10}$.

   A smooth complete intersection surface with Picard
   number one satisfies the assumptions of Theorem
   \ref{thm:irred-p3}. Thus by the Theorem of Noether the result
   applies in particular to general surfaces in $\PC^3$.

   \begin{corollary}
     Let $\Sigma\subset\PC^3$ be a smooth hypersurface of degree
     $n\geq 4$,
     let $H\subset\Sigma$ be a hyperplane
     section, and suppose that the Picard number  of
     $\Sigma$ is one.
     Let $d\geq n+6$ and let $\ks_1,\ldots,\ks_r$ be 
     topological (respectively analytical) singularity types.

     Suppose that 
     \begin{displaymath}
       \sum\limits_{i=1}^r \big(\tau^*(\ks_i)+2\big)^2
       <
       \tfrac{6\cdot\big(n^3-3n^2+8n-6\big)\cdot n^2}{\big(n^3-3n^2+10n-6\big)^2}
       \cdot (d+4-n)^2,
     \end{displaymath}

     Then 
     $V_{|dH|}^{irr}(\ks_1,\ldots,\ks_r)$ is non-empty
     and irreducible of the expected dimension.

     \begin{varthm-italic}[Addendum]
       If we, moreover, assume that $d\geq
       n\cdot\big(\tau^*(\ks_i)+1\big)$ for all $i=1,\ldots,r$, 
       then
       $V_{|dH|}^{irr}(\ks_1,\ldots,\ks_r)$ is also T-smooth.
     \hfill $\Box$
     \end{varthm-italic}
   \end{corollary}

   A general K3-surface has Picard number one and in this situation,
   by the Kodaira Vanishing Theorem
   $\Sigma$ also satisfies the assumption (ii) in Theorem \ref{thm:irred-p3}.

   \begin{corollary}
     Let $\Sigma$ be a smooth K3-surface with $\NS(\Sigma)=L\cdot\Z$
     with $L$ ample and set $n=L^2$.
     Let $d>0$, $D\sim_a dL$ and let $\ks_1,\ldots,\ks_r$ be 
     topological (respectively analytical) singularity types.

     Suppose that 
     \begin{displaymath}
       \sum\limits_{i=1}^r \big(\tau^*(\ks_i)+2\big)^2
       <
       \tfrac{54n^2+72n}{(11n+12)^2}\cdot d^2\cdot n.
     \end{displaymath}

     Then 
     $V_{|D|}^{irr}(\ks_1,\ldots,\ks_r)$ is irreducible and T-smooth,
     once it is non-empty. 

     \begin{varthm-italic}[Addendum]
       If  $d\geq 19$, then certainly
       $V_{|dH|}^{irr}(\ks_1,\ldots,\ks_r)$ is non-empty.
     \hfill $\Box$
     \end{varthm-italic}
   \end{corollary}


   \subsection{Products of Curves}\label{subsec:irreducibility:product-surfaces}

   If $\Sigma=C_1\times C_2$ is the product of two smooth projective
   curves, then for a general choice of $C_1$ and $C_2$ the
   N\'eron--Severi group will be generated by two fibres of the
   canonical projections, by abuse of notation also denoted by $C_1$
   and $C_2$. If both curves are elliptic, then ``general'' just means
   that the two curves are non-isogenous.

   \begin{theorem}\label{thm:irred-products-of-curves}
     Let $C_1$ and $C_2$ be two smooth projective curves of genera $g_1$
     and $g_2$  respectively with $g_1\geq g_2\geq 0$, such that for $\Sigma=C_1\times C_2$
     the N\'eron--Severi group is $\NS(\Sigma)=C_1\Z\oplus C_2\Z$.

     Let
     $\ks_1,\ldots,\ks_r$ be 
     topological or analytical singularity types, and
     let $D\in\Div(\Sigma)$ such that $D\sim_a aC_1+bC_2$ with
     \begin{displaymath}
       a\geq\left\{
       \begin{array}[m]{ll}
         \max\big\{2,\nu^*(\ks_i)\;\big|\;i=1,\ldots,r\big\}, &\text{
           if } g_2=0,\\
         2g_2-1, &\text{ else},
       \end{array}
       \right.
     \end{displaymath}
     and
     \begin{displaymath}
       b\geq\left\{
       \begin{array}[m]{ll}
         \max\big\{2,\nu^*(\ks_i)\;\big|\;i=1,\ldots,r\big\}, &\text{
           if } g_1=0,\\
         2g_1-1, &\text{ else}.
       \end{array}
       \right.
     \end{displaymath}

     Suppose that 
     \begin{equation}
       \label{eq:irred-products-of-curves:1}
       \sum\limits_{i=1}^r\big(\tau^*(\ks_i)+2\big)^2
       \;<\;\gamma\cdot(D-K_\Sigma)^2, 
     \end{equation} 
     where $\gamma$ may be taken from the following table with $\alpha=\frac{a-2g_2+2}{b-2g_1+2}>0$.
     \as{1.3}
     \begin{center}
       \tablefirsthead{\hline
         \multicolumn{1}{|c}{$g_1$}&\multicolumn{1}{|c|}{$g_2$}&$\gamma$\\\hline\hline}
       \tablehead{\hline
         \multicolumn{1}{|c}{$g_1$}&\multicolumn{1}{|c|}{$g_2$}&$\gamma$\\\hline\hline}
       \tabletail{\hline}
       \tablelasttail{\\\hline}
       \begin{supertabular}{|r|r|c|}         
         $0$&$0$&$\frac{1}{24} $\\
         $1$&$0$&$\frac{1}{\max\{32,2\alpha\}}$\\
         $\geq 2$&$0$&
         $\frac{1}{\max\{24+16g_1,4g_1\alpha\}} $\\
         $1$&$1$&
         $\frac{1}{\max\big\{32,2\alpha,\tfrac{2}{\alpha}\big\}} $\\
         $\geq 2$&$\geq 1$&
         $\frac{1}{\max\left\{24+16g_1+16g_2,4g_1\alpha,\tfrac{4g_2}{\alpha}\right\}}$
       \end{supertabular}
     \end{center}
     \as{1}

     Then 
     $V_{|D|}^{irr}(\ks_1,\ldots,\ks_r)$ is empty
     or  irreducible of the expected dimension.
%
     \hfill $\Box$
   \end{theorem}

   Only in the case $\Sigma\cong\PC^1\times\PC^1$ we get a
   constant $\gamma$ which does not depend on the chosen divisor
   $D$, while in the remaining cases  the ratio of $a$ and $b$ is
   involved in $\gamma$. This means that an asymptotical behaviour can
   only be examined if the ratio remains unchanged.


   \subsection{Geometrically Ruled Surfaces}\label{subsec:irreducibility:ruled-surfaces}

   Let $\pi:\Sigma=\P(\ke)\rightarrow C$ be a geometrically ruled
   surface with normalised bundle $\ke$ (in the 
   sense of \cite{Har77} V.2.8.1). The N\'eron--Severi group of
   $\Sigma$ is $\NS(\Sigma) = C_0\Z\oplus F\Z$ 
   with intersection matrix $\left(\begin{smallmatrix}-e & 1 \\ 1 & 0\end{smallmatrix}\right)$
   where $F\cong\PC^1$ is a fibre of $\pi$, $C_0$ a section of $\pi$
   with $\ko_\Sigma(C_0)\cong\ko_{\P(\ke)}(1)$, $g=g(C)$ the genus of
   $C$,  $\mathfrak{e}=\Lambda^2\ke$ and
   $e=-\deg(\mathfrak{e})\geq -g$. 
   For the canonical divisor we have $K_\Sigma \sim_a -2C_0+ (2g-2-e)\cdot F$.

   \begin{theorem}\label{thm:irred-ruled-surfaces}
     Let $\pi:\Sigma\rightarrow C$ be a geometrically ruled surface
     with  $e\leq 0$.
     Let
     $\ks_1,\ldots,\ks_r$ be 
     topological or analytical singularity types,
     and let $D\in\Div(\Sigma)$ such that $D\sim_a aC_0+bF$ with
     $a\geq\max\big\{2,\nu^*(\ks_i)\;\big|\;i=1,\ldots,r\big\}$, 
     and,
     \begin{displaymath}
       b>\left\{
         \begin{array}{ll}
           \max\big\{1,\nu^*(\ks_i)-1\;\big|\;i=1,\ldots,r\big\}, &
           \mbox{ if } g=0,\\
           2g-2+\frac{ae}{2}, & \mbox{ if } g>0.
         \end{array}
       \right.
     \end{displaymath}

     Suppose that 
     \begin{equation}
       \label{eq:irred-ruled-surfaces:1}
       \sum\limits_{i=1}^r\big(\tau^*(\ks_i)+2\big)^2
       \;<\;\gamma\cdot(D-K_\Sigma)^2, 
     \end{equation} 
     where $\gamma$ may be taken from the following table with
     $\alpha=\frac{a+2}{b+2-2g-\tfrac{ae}{2}}>0$. 
     \as{1.3}
     \begin{center}
       \begin{tabular}{|r|r|c|}
         \hline
         \multicolumn{1}{|c}{$g$}&\multicolumn{1}{|c|}{$e$}&$\gamma$\\\hline\hline
          $0$ & $0$ &  $\tfrac{1}{24}$ \\
          $1$ & $0$ & $\tfrac{1}{\max\{24,2\alpha\}}$ \\
          $1$ & $-1$ & $\tfrac{1}{\max\Big\{\min\big\{30+\tfrac{16}{\alpha}+4\alpha,40+9\alpha\big\},
           \tfrac{13}{2}\alpha\Big\}}$\\ 
          $\geq 2 $&$ 0 $&$  \tfrac{1}{\max\{24+16g, 4g\alpha\}}$\\ 
          $\geq 2 $&$ <0 $& 
          $ \tfrac{1}{\max\Big\{\min\big\{24+16g-9e\alpha,18+16g-9e\alpha
           -\tfrac{16}{e\alpha}\big\},4g\alpha-9e\alpha\Big\}}$ \\\hline
       \end{tabular}
     \end{center}
     \as{1}

     Then 
     $V_{|D|}^{irr}(\ks_1,\ldots,\ks_r)$ is empty
     or  irreducible of the expected dimension.
     \hfill $\Box$
   \end{theorem}

   Once more, only in the case $g=0$, i.\ e.\ when $\Sigma\cong\PC^1\times\PC^1$, we
   are in the lucky situation that the 
   constant $\gamma$ does not at all depend on the chosen divisor
   $D$, whereas in the case $g\geq 1$ the ratio of $a$ and $b$ is
   involved in $\gamma$. This means that an asymptotical behaviour can
   only be examined if the ratio remains unchanged.

   If $\Sigma$ is a product $C\times\PC^1$ the constant $\gamma$ here is
   the same as in Section 
   \ref{subsec:irreducibility:product-surfaces}. 


   \section{The Proofs}\label{sec:technical-lemmata}

   Our approach to the problem proceeds along the lines of an unpublished result of
   Greuel, Lossen and Shustin (cf.~\cite{GLS98d}), which is based on
   ideas of Chiantini and Ciliberto (cf.~\cite{CC99}). It is a slight
   modification of the proof given in \cite{Kei03}.
   
   We tackle the problem in three steps:\label{page:fix}
   \begin{varthm-roman}[Step 1]
     By \cite{Kei03} Theorem 3.1 we know  that the open subvariety $V^{irr,reg}$ of curves in
     $V^{irr}$ with 
     $h^1\big(\Sigma,\kj_{X(C)/\Sigma}(D)\big)=0$ is always irreducible,
     and hence so is its 
     closure in $V^{irr}$. 
   \end{varthm-roman}
   \begin{varthm-roman}[Step 2]
     We find conditions which ensure that the open subvariety
     $V^{irr,fix}$ of curves in $V^{irr}$ with 
     $h^1\big(\Sigma,\kj_{X^*_{fix}(C)/\Sigma}(D)\big)=0$ is dense in
     $V^{irr}$. 
   \end{varthm-roman}
   \begin{varthm-roman}[Step 3]
     And finally, we combine these conditions with conditions which
     guarantee that $V^{irr,reg}$ is dense in $V^{irr,fix}$ by showing
     that they share some open dense subset $V^{gen}_U$ of curves with
     singularities in very general position.      
 
     More precisely, taking Lemma \ref{lem:irred-F} into account, we deduce from
     Lemma \ref{lem:irred-G}  
     conditions which ensure that there exists a very 
     general subset $U\subset\Sigma^r$  such
     that the family
     $V^{gen}_U=V^{gen}_{|D|,U}(\ks_1,\ldots,\ks_r)$, as defined
     there, satisfies
     \begin{enumerate}
     \item $V^{gen}_U$ is dense in $V^{irr,fix}$, and
     \item $V^{gen}_U\subseteq V^{irr,reg}$.
     \end{enumerate}

     But then $V^{irr,reg}$ is dense in $V^{irr}$ and $V^{irr}$ is
     irreducible by Step 1.
     \hfill$\Box$
  \end{varthm-roman}

   The difficult part is Step 2. 
   For this one we consider
   the restriction of the 
   morphism (cf.\ Subsection \ref{subsec:families})
   \begin{displaymath}
     \Phi:V^{irr}\rightarrow \Sym^r(\Sigma)=:\kb
   \end{displaymath}
   to an irreducible component $V^*$ of $V^{irr}$ not contained in the
   closure 
   $\overline{V^{irr,fix}}$ in $V^{irr}$. Knowing, that the dimension of $V^*$ is
   at least the expected dimension $\dim\big(V^{irr,fix}\big)$ we
   deduce that the codimension of $\Phi\big(V^*\big)$ in $\kb$ is 
   at most $h^1\big(\Sigma,\kj_{X^*_{fix}(C)/\Sigma}(D)\big)$, where $C\in
   V^*$ (cf.~Lemma \ref{lem:irred-E-fix}). It thus suffices to find
   conditions which contradict this 
   inequality, that is, we have to get our hands on
   $\codim_\kb\big(\Phi(V^*)\big)$. This is achieved  by applying the results of
   \cite{Kei03} Lemma 4.1 to Lemma 4.6 to the zero-dimensional scheme
   $X_0=X^*_{fix}(C)$. 

   These considerations lead to the following proofs.

   \begin{proof}[Proof of Theorem \ref{thm:irred-p3}]
     We may assume that
     $V^{irr}=V_{|D|}^{irr}(\ks_1,\ldots,\ks_r)$
     is non-empty.
     As indicated above it suffices to show that:
     \begin{varthm-roman}[Step 2] $V^{irr}=\overline{V^{irr,fix}}$, where 
       $V^{irr,fix}=V_{|D|}^{irr,fix}(\ks_1,\ldots,\ks_r)$, and
     \end{varthm-roman}
     \begin{varthm-roman}[Step 3] the conditions of Lemma \ref{lem:irred-G} are fulfilled.
     \end{varthm-roman}

     For Step 3 we note that
     $\nu^*(\ks_i)\leq\tau^*(\ks_i)$.
     Thus \eqref{eq:irred-p3:2} implies that
     \begin{displaymath}
       \sum_{i=1}^r  \big(\nu^*(\ks_i)+2\big)^2
       \leq
       \sum_{i=1}^r  \big(\tau^*(\ks_i)+2\big)^2
       \leq
       \gamma\cdot (D-K_\Sigma)^2
       \leq
       \tfrac{1}{2}\cdot (D-K_\Sigma)^2,
     \end{displaymath}
     which gives the first condition in Lemma \ref{lem:irred-G}. Since
     a surface with Picard number one has no curves of
     selfintersection zero, the second condition in Lemma
     \ref{lem:irred-G} is void, while the last condition is satisfied
     by \eqref{eq:irred-p3:0}.

     It remains to show Step 2, i.~e.~$V^{irr}=\overline{V^{irr,fix}}$.
     Suppose the contrary, that is, there is an irreducible curve
     $C_0\in V^{irr}\setminus\overline{V^{irr,fix}}$, in particular
     $h^1\big(\Sigma,\kj_{X_0/\Sigma}(D)\big)>0$ for $X_0=X^*_{fix}(C_0)$.
     Since
     $\deg(X_0)
     =\sum_{i=1}^r\big(\tau^*(\ks_i)+2\big)$ and 
     $\sum_{z\in\Sigma}\big(\deg(X_{0,z})\big)^2=
     \sum_{i=1}^r\big(\tau^*(\ks_i)+2\big)^2$ the assumptions (0)-(3) of
     \cite{Kei03} Lemma 4.1 and (4) of \cite{Kei03} Lemma 4.3 are
     fulfilled. Thus \cite{Kei03} Lemma 4.3 implies that $C_0$
     satisfies Condition \eqref{eq:irred-E-fix:1} in Lemma
     \ref{lem:irred-E-fix} below, which it cannot
     satisfy by the same Lemma. Thus we have derived a
     contradiction.
   \end{proof}

   \begin{proof}[Proof of Theorem \ref{thm:irred-products-of-curves}]
     The assumptions on $a$ and $b$ ensure that $D-K_\Sigma$ is big
     and nef and that $D+K_\Sigma$ is nef. Thus,
     once we know that \eqref{eq:irred-products-of-curves:1} implies Condition (3)
     in \cite{Kei03} Lemma 4.1 we can do the same proof as in Theorem
     \ref{thm:irred-p3}, just replacing \cite{Kei03} Lemma 4.3 by
     \cite{Kei03} Lemma 4.4.

     For Condition (3) we note that
     \begin{displaymath}
         \sum\limits_{i=1}^r\deg\big(X^*_{fix}(\ks_i)\big)
         \leq\sum\limits_{i=1}^r\big(\tau^*(\ks_i)+2\big)^2 
         \leq\tfrac{1}{24}\cdot(D-K_\Sigma)^2
         <\tfrac{1}{4}\cdot(D-K_\Sigma)^2.
     \end{displaymath}
   \end{proof}

   \begin{proof}[Proof of Theorem \ref{thm:irred-ruled-surfaces}]
     The proof is identical to that of Theorem
     \ref{thm:irred-products-of-curves}, just replacing 
     \cite{Kei03} Lemma 4.4 by \cite{Kei03} Lemma 4.6.
   \end{proof}


   \subsection{Some Technical Lemmata}

   Can have applied \cite{Kei03} Lemma 4.1 to Lemma 4.6 to the
   zero-dimensional scheme $X_0=X^*_{fix}(C)$, for a curve
   $C\in V^{irr}\setminus\overline{V^{irr,fix}}$, in order to find with
   the aid of Bogomolov instability curves $\Delta_i$ and subschemes
   $X_i^0\subseteq X_i$, where $X_i=X_{i-1}:\Delta_i$, such that for
   $X_S=\bigcup_{i=1}^mX_i^0$ 
   \begin{displaymath}
     h^1\big(\Sigma,\kj_{X_0}(D)\big)+
     \sum_{i=1}^m
     \Big(h^0\big(\Sigma,\ko_\Sigma(\Delta_i)\big)-1\Big) < \# X_S.
   \end{displaymath}
   And we are now going to show that this simply is not possible.

   \begin{lemma}\label{lem:irred-E-fix}
     Let $D\in\Div(\Sigma)$, $\ks_1,\ldots,\ks_r$ be pairwise distinct
     topological (respectively analytical) singularity types. Suppose that
     $V_{|D|}^{irr,fix}(\ks_1,\ldots,\ks_r)$ is non-empty.

     Then there exists no curve\footnote{For a subset $U\subseteq V$
       of a topological space $V$ we denote by $\overline{U}$ the
       closure of $U$ in $V$.} 
     $C\in V_{|D|}^{irr}(\ks_1,\ldots,\ks_r)
     \setminus
     \overline{V_{|D|}^{irr,fix}(\ks_1,\ldots,\ks_r)}$
     such that for the zero-dimensional scheme $X_0=X^*_{fix}(C)$
     there exist curves 
     $\Delta_1,\ldots,\Delta_m\subset\Sigma$ and
     zero-dimensional locally complete intersections $X_i^0\subseteq
     X_{i-1}$ for $i=1,\ldots,m$, where $X_i=X_{i-1}:\Delta_i$ for
     $i=1,\ldots,m$ 
     such that $X_S=\bigcup_{i=1}^mX_i^0$ satisfies
     \begin{equation}
       \label{eq:irred-E-fix:1}
       h^1\big(\Sigma,\kj_{X_0}(D)\big)+
       \sum_{i=1}^m
       \Big(h^0\big(\Sigma,\ko_\Sigma(\Delta_i)\big)-1\Big) < \# X_S.
     \end{equation}
   \end{lemma}

   \begin{proof}
     Throughout the proof we use the notation
     $V^{irr}=V_{|D|}^{irr}(\ks_1,\ldots,\ks_r)$
     and
     $V^{irr,fix}=V_{|D|}^{irr,fix}(\ks_1,\ldots,\ks_r)$.
     
     Suppose there exists a curve $C\in V^{irr}\setminus \overline{V^{irr,fix}}$
     satisfying the assumption of the Lemma, and let $V^*$ be the
     irreducible component of $V^{irr}$
     containing $C$. Moreover, let $C_0\in V^{irr,fix}$.

     We consider in the following the morphism
     from Subsection \ref{subsec:families}
     \begin{displaymath}
       \Phi=\Phi_{|D|}(\ks_1,\ldots,\ks_r)
       :V_{|D|}(\ks_1,\ldots,\ks_r)
       \rightarrow
       \Sym^r(\Sigma)=:\kb.
     \end{displaymath}

     \begin{varthm-roman}[Step 1]
       $h^0\big(\kj_{X^*_{fix}(C_0)/\Sigma}(D)\big)
       =
       h^0\big(\kj_{X^*_{fix}(C)/\Sigma}(D)\big)-
       h^1\big(\kj_{X^*_{fix}(C)/\Sigma}(D)\big).$
     \end{varthm-roman}
     By the choice of $C_0$ we have
     \begin{displaymath}
       0=H^1\big(\Sigma,\kj_{X^*_{fix}(C_0)/\Sigma}(D)\big)\rightarrow
       H^1(\Sigma,\ko_\Sigma(D)\big) \rightarrow
       H^1(\Sigma,\ko_{X^*_{fix}(C_0)}(D)\big)=0,
     \end{displaymath}
     and thus $D$ is non-special,
     i.~e.~$h^1(\Sigma,\ko_\Sigma(D)\big)=0$.
     But then
     \begin{displaymath}
       \begin{array}{rcl}
         h^0\big(\Sigma,\kj_{X^*_{fix}(C_0)/\Sigma}(D)\big)
         &=&h^0\big(\Sigma,\ko_\Sigma(D)\big)-\deg\big(X^*_{fix}(C_0)\big)\\
         &=&h^0\big(\Sigma,\ko_\Sigma(D)\big)-\deg\big(X^*_{fix}(C)\big)\\
         &=&
         h^0\big(\Sigma,\kj_{X^*_{fix}(C)/\Sigma}(D)\big)-
         h^1\big(\Sigma,\kj_{X^*_{fix}(C)/\Sigma}(D)\big).
       \end{array}
     \end{displaymath}

     \begin{varthm-roman}[Step 2]       
       $h^1\big(\Sigma,\kj_{X^*_{fix}(C)}(D)\big)\geq\codim_\kb\Big(\Phi\big(V^*\big)\Big)$.
     \end{varthm-roman}
     Suppose the contrary, that is
     $\dim\Big(\Phi\big(V^*\big)\Big)<\dim(\kb)-h^1\big(\Sigma,\kj_{X^*_{fix}(C)/\Sigma}(D)\big)$.
     The vanishing of
     $h^1\big(\Sigma,\kj_{X^*_{fix}(C_0)/\Sigma}(D)\big)$ implies that
     $V^{irr}$ is smooth of the expected dimension
     $\dim\big(V^{irr,fix}\big)$ at
     $C_0$. Therefore, and
     by Step 1 and  Equation
     \eqref{eq:psi-fix}  we have
     \begin{displaymath}
       \begin{array}{rcl}
         \dim\big(V^*\big)&\leq& \dim\Big(\Phi\big(V^*\big)\Big)+
         \dim\Big(\Phi^{-1}\big(\Phi(C)\big)\Big)\\
         &<&\dim(\kb)-h^1\big(\Sigma,\kj_{X^*_{fix}(C)/\Sigma}(D)\big)+
         h^0\big(\Sigma,\kj_{X^*_{fix}(C)/\Sigma}(D)\big)-1\\
         &=&
         2r+
         h^0\big(\Sigma,\kj_{X^*_{fix}(C_0)/\Sigma}(D)\big)-1=\dim\big(V^{irr,fix}\big).         
       \end{array}
     \end{displaymath}
     However, any irreducible component of $V^{irr}$ has at least the
     expected dimension $\dim\big(V^{irr,fix}\big)$, which gives a
     contradiction. 

     \begin{varthm-roman}[Step 3]
       $\codim_\kb\Big(\Phi\big(V^*\big)\Big)\geq\#X_S-\sum_{i=1}^m\dim|\Delta_i|_l$.
     \end{varthm-roman}
     The existence of the subschemes $X_i^0\subseteq X^*_{fix}(C)\cap
     \Delta_i$ imposes at least $\# X_i^0-\dim|\Delta_i|_l$ conditions
     on $X^*_{fix}(C)$ and increases thus the codimension of
     $\Phi\big(V^*\big)$ by the same number. 

     \begin{varthm-roman}[Step 4]
       Derive a contradiction.
     \end{varthm-roman}
     Collecting the results we derive the following contradiction:
     \begin{multline*}
         h^1\big(\Sigma,\kj_{X^*_{fix}(C)}(D)\big)\geq_\expl{Step 2}
         \codim_\kb\Big(\Phi\big(V^*\big)\Big)\\
         \geq_\expl{Step 3}
         \#X_S-\mbox{$\sum_{i=1}^m$}\dim|\Delta_i|_l
         >_\expl{\eqref{eq:irred-E-fix:1}}h^1\big(\Sigma,\kj_{X^*_{fix}(C)}(D)\big).
     \end{multline*}
   \end{proof}

   The next two lemmata provide conditions which ensure that
   $V^{irr,reg}$ and $V^{irr,fix}$ share some dense subset
   $V^{gen}_U$, and thus that $V^{irr,reg}$ is dense in $V^{irr,fix}$.

   \begin{lemma}\label{lem:irred-F}
     Let $\ks_1,\ldots,\ks_r$ be topological (respectively analytical) singularity
     types, let $D\in\Div(\Sigma)$ and let $V^{irr}=
     V_{|D|}^{irr}(\ks_1,\ldots,\ks_r)$.

     There exists a very general subset $U\subset\Sigma^r$  
     such that\footnote{Here $\sim$ means either
       topological equivalence $\sim_t$
       or contact equivalence $\sim_c$. -- By a very general subset of
       $\Sigma^r$ we mean the complement of at most countably many
       closed subvarieties.}
     $V^{gen}_U=V^{gen}_{|D|,U}(\ks_1,\ldots,\ks_r) =
     \big\{C\in V^{irr}\;\big|\; \underline{z}\in U, (C,z_i)\sim S_i,i=1,\ldots,r\}$ 
     is dense in $V^{irr,fix}_{|D|}(\ks_1,\ldots,\ks_r)$.
   \end{lemma}
   \begin{proof}
     This follows from the remark in Subsection \ref{subsec:families}. 
   \end{proof}

   \begin{lemma}\label{lem:irred-G}
     With the notation of Lemma \ref{lem:irred-F} we assume that
     \begin{enumerate}
     \item $(D-K_\Sigma)^2\geq
       2\cdot\sum_{i=1}^k\big(\nu^*(\ks_i)+2\big)^2$,
     \item
       $(D-K_\Sigma).B>\max\big\{\nu^*(\ks_i)+1\;\big|\;i=1,\ldots,r\big\}$\; for any
          irreducible curve $B$ with $B^2=0$ and $\dim|B|_a>0$, and 
     \item $D-K_\Sigma$ is nef.
     \end{enumerate}
     Then there exists a very general subset $U\subset\Sigma^r$  such that 
     $V^{gen}_U\subseteq
     V^{irr,reg}_{|D|}(\ks_1,\ldots,\ks_r)$.
   \end{lemma}
   \begin{proof}
     By \cite{KT02} Theorem 2.1 we know that there is a very
     general subset $U\subset\Sigma^r$ such that for
     $\underline{z}\in U$ and
     $\underline{\nu}=\big(\nu^*(\ks_1)+1,\ldots,\nu^*(\ks_r)+1\big)$  we have
     \begin{displaymath}
       h^1\big(\Sigma,\kj_{X(\underline{\nu};\underline{z})/\Sigma}(D)\big)=0.
     \end{displaymath}
     However, if $C\in V^{irr}$ and $\underline{z}\in U$ with
     $(C,z_i)\sim \ks_i$, then by the definition of $\nu^*(\ks_i)$ we
     have 
     \begin{displaymath}
       \kj_{X(\underline{\nu};\underline{z})/\Sigma}\hookrightarrow 
       \kj_{X(C)/\Sigma},
     \end{displaymath}
     and hence 
     the vanishing of
     $H^1\big(\Sigma,\kj_{X(\underline{\nu};\underline{z})/\Sigma}(D)\big)$
     implies 
       $h^1\big(\Sigma,\kj_{X(C)/\Sigma}(D)\big)=0$,
     i.~e.~$C\in V^{irr,reg}$.
   \end{proof}

   \bibliographystyle{amsalpha-tom}
   \bibliography{bibliographie}

\providecommand{\bysame}{\leavevmode\hbox to3em{\hrulefill}\thinspace}
\begin{thebibliography}{ChC99}

\bibitem[ChC99]{CC99}
Luca Chiantini and Ciro Ciliberto, \emph{On the {Severi} variety of surfaces in
  {$\mathbbm P_{\mathbbm C}^3$}}, J.\ Algebraic Geom. \textbf{8} (1999),
  67--83.

\bibitem[GLS97]{GLS97}
{Gert-Martin} Greuel, Christoph Lossen, and Eugenii Shustin, \emph{New
  asymptotics in the geometry of equisingular families of curves}, Internat.\
  Math.\ Res.\ Notices \textbf{13} (1997), 595--611.

\bibitem[GLS98]{GLS98d}
{Gert-Martin} Greuel, Christoph Lossen, and Eugenii Shustin, \emph{On the
  irreducibility of families of curves}, Unpulished Manuscript, 1998.

\bibitem[GLS00]{GLS00}
{Gert-Martin} Greuel, Christoph Lossen, and Eugenii Shustin, \emph{Castelnuovo
  function, zero-dimensional schemes, and singular plane curves}, J.\ Algebraic
  Geom. \textbf{9} (2000), no.~4, 663--710.

\bibitem[Har77]{Har77}
Robin Hartshorne, \emph{Algebraic geometry}, Springer, 1977.

\bibitem[Kei02]{Kei03}
Thomas Keilen, \emph{Irreducibility of equisingular families of curves},
  accepted by Trans.\ Amer.\ Math.\ Soc., http:// \!\!www. \!\!mathematik.
  \!\!uni-kl. \!\!de/ \!\!\textasciitilde keilen/ \!\!download/ \!\!Keilen001/
  \!\!Keilen001.ps.gz, 2002.

\bibitem[KeT02]{KT02}
Thomas Keilen and Ilya Tyomkin, \emph{Existence of curves with prescribed
  singularities}, Trans.\ Amer.\ Math.\ Soc. \textbf{354} (2002), no.~5,
  1837--1860, http:// \!\!www. \!\!mathematik. \!\!uni-kl. \!\!de/
  \!\!\textasciitilde keilen/ \!\!download/ \!\!KeilenTyomkin001/
  \!\!KeilenTyomkin001.ps.gz.

\bibitem[Los98]{Los98}
Christoph Lossen, \emph{The geometry of equisingular and equianalytic families
  of curves on a surface}, Phd thesis, FB Mathematik, Universit\"at
  Kaiserslautern, Aug. 1998, http:// \!\!www. \!\!mathematik. \!\!uni-kl.
  \!\!de/ \!\!\textasciitilde lossen/ \!\!download/ \!\!Lossen002/
  \!\!Lossen002.ps.gz.

\end{thebibliography}

\Addresses

\end{document}